\documentclass[12pt]{article}
\usepackage{amssymb,amsfonts,amsmath, psfrag,eepic,colordvi,graphicx,epsfig,ytableau}
\usepackage{amssymb,latexsym,graphics}
\usepackage{MnSymbol}
\usepackage[enableskew]{youngtab}
\usepackage{float}
\usepackage{tikz}
\topskip  
\newcommand{\rmnum}[1]{\romannumeral #1}

 \numberwithin{equation}{section}
\newtheorem{theorem}{Theorem}[section]

\newtheorem{conjecture}[theorem]{Conjecture}
\newtheorem{remark}[theorem]{Remark}
\newtheorem{lemma}[theorem]{Lemma}

\newtheorem{observation}[theorem]{Observation}

\begin{document}
\parskip 6pt

\pagenumbering{arabic}
\def\sof{\hfill\rule{2mm}{2mm}}
\def\ls{\leq}
\def\gs{\geq}
\def\SS{\mathcal S}
\def\qq{{\bold q}}
\def\MM{\mathcal M}
\def\TT{\mathcal T}
\def\EE{\mathcal E}
\def\lsp{\mbox{lsp}}
\def\rsp{\mbox{rsp}}
\def\pf{\noindent {\it Proof.} }
\def\mp{\mbox{pyramid}}
\def\mb{\mbox{block}}
\def\mc{\mbox{cross}}
\def\qed{\hfill \rule{4pt}{7pt}}
\def\block{\hfill \rule{5pt}{5pt}}

\begin{center}
{\Large\bf  Partial $\gamma$-Positivity for   Quasi-Stirling Permutations of  Multisets}
\end{center}

\begin{center}
{\small Sherry H.F. Yan,   Yunwei Huang,  Lihong Yang }

 Department of Mathematics\\
Zhejiang Normal University\\
 Jinhua 321004, P.R. China

 hfy@zjnu.cn

\end{center}

\noindent {\bf Abstract.}
We prove that the enumerative polynomials of quasi-Stirling permutations of multisets with respect to the statistics of plateaux, descents and ascents are partial $\gamma$-positive, thereby confirming a recent conjecture posed by Lin, Ma and Zhang. This is accomplished by proving the partial $\gamma$-positivity of the enumerative polynomials of certain ordered labeled trees, which are in bijection with quasi-Stirling permutations of multisets.  As an application, we  provide  an alternative proof of the  partial $\gamma$-positivity of the enumerative polynomials  on    Stirling permutations of multisets.




\section{Introduction}
Gamma-positive polynomials arise frequently in combinatorics; the reader is referred to the survey of
Athanasiadis \cite{Ath} for more information.  A univariate
polynomial $f_n(x)$ is said to be {\em $\gamma$-positive} if it can be expanded as
$$
f_n(x)=\sum_{k=0}^{\lfloor {n\over 2}\rfloor} \gamma_k x^k(1+x)^{n-2k}
$$
with $\gamma_k\geq 0$. A bivariate polynomial $h_n(x,y)$ is said to be {\em homogeneous $\gamma$-positive}, if  $h_n(x, y)$ can be expressed as
$$
h_n(x,y)=\sum_{k=0}^{\lfloor {n\over 2}\rfloor} \gamma_k (xy)^k(x+y)^{n-2k}
$$
with $\gamma_k\geq 0$.
 Foata and Sch\"uzenberger \cite{Fo} proved that the bivariate Eulerian polynomial $A_n(x,y)=\sum_{\pi\in \mathcal{S}_n} x^{asc(\pi)}y^{des(\pi)}$ is homogeneous $\gamma$-positive, where $\mathcal{S}_n$ is the set of permutations of $[n]=\{1,2,\ldots, n\}$.   A trivariate polynomial
$p(x, y, z)=\sum_{i}s_{i}(x,y)z^i $
is said to be {partial $\gamma$-positive} if every $s_i(x, y)$ is homogeneous
$\gamma$-positive.

Let $\mathcal{M}=\{1^{k_1}, 2^{k_2}, \ldots, n^{k_n}\}$ be a multiset where $k_i$ is the number of occurrences of $i$ in  $\mathcal{M}$ and $k_i\geq 1$.
  A permutation $\pi=\pi_1\pi_2\ldots \pi_n$ of a multiset $\mathcal{M}$ is said to be  a {\em Stirling } permutation if   $i<j<k$ and $\pi_i=\pi_k$, then $\pi_j>\pi_i$. For a multiset $\mathcal{M}$, let $ \mathcal{Q}_{\mathcal{M}}$ denote  the set of Stirling permutations of $\mathcal{M}$. For example, if $\mathcal{M}=\{1^2, 2^2\}$, we have
$
 \mathcal{Q}_{\mathcal{M}} =\{1221, 2211, 1122\}.
$
  Stirling permutations were originally introduced by Gessel and Stanley \cite{Gessel} in the case of the multiset $\mathcal{M}=\{1^{2}, 2^{2}, \ldots, n^{2}\}$. There is an extensive literature on Stirling permutations, see   \cite{Brenti1, Brenti2, Dzhumadil'daev, Liu, Park, Park2} and references therein.
In   analogy to Stirling permutations, Archer et al. \cite{Archer} introduced   {\em quasi-Stirling} permutations. A permutation $\pi=\pi_1\pi_2\ldots \pi_n$ of a multiset is said to be a {\em quasi-Stirling
} permutation  if  there does not exist four indices $i<j<k<\ell$ such that $\pi_i=\pi_k$ and $\pi_j=\pi_\ell$. For a multiset $\mathcal{M}$, denote by   $\overline{\mathcal{Q}}_{\mathcal{M}}$ the set of quasi-Stirling permutations of $\mathcal{M}$. For example, if $\mathcal{M}=\{1^2, 2^2\}$, we have
$
\overline{\mathcal{Q}}_{\mathcal{M}} =\{1221, 2112, 1122, 2211\}.
$
Clearly, we have $\mathcal{Q}_{\mathcal{M}}\subseteq \overline{\mathcal{Q}}_{\mathcal{M}} $ for a fixed multiset $\mathcal{M}$.

For a permutation $\pi=\pi_1\pi_2\ldots \pi_n$,   an index $i$,  $0\leq i\leq n$, is called an  {\em ascent} (resp. {\em a descent},  {\em a plateau}) of $\pi$ if $\pi_{i}<\pi_{i+1}$ (resp. $\pi_i>\pi_{i+1}$, $\pi_i=\pi_{i+1}$) with the convention $\pi_0=\pi_{n+1}=0$.  Let $asc(\pi)$ (resp. $des(\pi)$, $plat(\pi)$) denote the number of   ascents (resp. descents, plateaux) of $\pi$.
An index   $i$, $1\leq i\leq n$, is called a {\em double descent} of $\pi$ if $\pi_{i-1}>\pi_{i}>\pi_{i+1}$ with the convention $\pi_0=\pi_{n+1}=0$.  Let   $ddes(\pi)$  denote   the number  of   double descents of $\pi$.
In \cite{Bona}, Bona proved that   the three
statistics $asc$, $plat$ and $des$ are equidistributed over Stirling permutations on $\mathcal{M}=\{1^2, 2^2, \ldots, n^2\}$. In \cite{Elizalde},    Elizalde  investigated the distribution of the  statistics of plateaux, descents and ascents  on  quasi-Stirling permutations of the multiset $\mathcal{M}=\{1^k, 2^k, \ldots, n^k\}$.

For a multiset $\mathcal{M}$, we
define $$ \overline{Q}_{\mathcal{M}}(x,y,z)=\sum_{\pi\in \overline{\mathcal{Q}}_{\mathcal{M}}}x^{asc(\pi)}y^{des(\pi)}z^{plat(\pi)},$$
and
$$ Q_{\mathcal{M}}(x,y,z)=\sum_{\pi\in \mathcal{Q}_{\mathcal{M}}}x^{asc(\pi)}y^{des(\pi)}z^{plat(\pi)}.$$

In \cite{Lin}, Lin, Ma and Zhang  proved that the polynomial  $Q_{\mathcal{M}}(x,y,z)$ is partial $\gamma$-positive via the machine of context-free
grammars and a group action on Stirling multipermutations. Specializing their result  to the  so-called Jacobi-Stirling permutations introduced in \cite{Gesse2} confirms  a recent partial $\gamma$-positivity conjecture due to Ma, Ma and Yeh \cite{Ma}. Lin, Ma and Zhang \cite{Lin} also posed the following conjecture.

\begin{conjecture}\label{con1}
For a multiset $\mathcal{M}$,
the polynomial
$\overline{Q}_{\mathcal{M}}(x,y,z)$ is  partial $\gamma$-positive.
\end{conjecture}

In this paper,  we confirm Conjecture \ref{con1} by proving  the partial $\gamma$-positivity of the enumerative polynomials of certain ordered labeled trees, which are shown  to be in bijection with quasi-Stirling permutations in \cite{Yanzhu}.  As an application, we  provide  an alternative proof of the  partial $\gamma$-positivity of  the polynomial $Q_{\mathcal{M}}(x,y,z)$.

The rest of this paper is organized as follows. In Section 2, we prove the partial $\gamma$-positivity of the enumerative polynomials of certain ordered labeled trees by introducing a group action  in the sprit of the  Foata-Strehl
action  on permutations \cite{Fo}.   Section 3 is devoted to the proofs of the    partial $\gamma$-positivity of  the polynomial $\overline{Q}_{\mathcal{M}}(x,y,z)$ and $Q_{\mathcal{M}}(x,y,z)$.

\section{Partial $\gamma$-positivity for Trees}

The objective of this section is  to prove the partial $\gamma$-positivity of the enumerative polynomials of certain ordered labeled trees, which were  shown to be in bijection with quasi-Stirling permutations in \cite{Yanzhu}.
Recall that an  {\em ordered} tree is a tree with one designated vertex, which is called the root, and the subtrees of
each vertex are linearly ordered.  In an ordered tree $T$,   the {\em level} of a vertex $v$ in $T$  is defined to be the length of the unique path from the root to $v$.  A vertex $v$ is said to be {\em   odd  } (resp. {\em even}) if the level of $v$ is odd (resp. even).

\begin{figure}[h]
\centerline{\includegraphics[width=5cm]{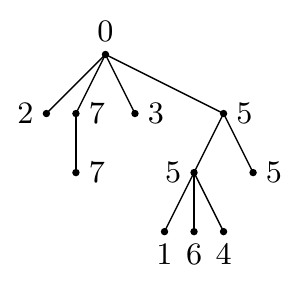}}
    \caption{A tree $\mathcal{T}\in \mathcal{T}_\mathcal{M}$ with $\mathcal{M}=\{1,2,3,4,5^3,6,7^2\}$}\label{T}
\end{figure}

  Yan and Zhu \cite{Yanzhu}  introduced   a new
class of ordered   labeled trees $T$   verifying the following properties:
\begin{itemize}
\item[{\upshape (\rmnum{1})}] the vertices are labeled by the elements of the multiset  $\{0\}\cup \mathcal{M}$, where $\mathcal{M}=\{1^{k_1}, 2^{k_2}, \ldots, n^{k_n}\}$ with $k_i\geq 1$;
    \item[{\upshape (\rmnum{2})}] the root is labeled by $0$;
    \item[{\upshape (\rmnum{3})}] for  an odd vertex  $v$,  if $v$ is labeled by $i$, then $v$ has  exactly  $k_i-1$ children  and the children of $v$ have the same label as that of $v$
  in $T$.
\end{itemize}
Let $\mathcal{T}_{\mathcal{M}}$ denote the set of  such ordered    labeled trees. For example,  a tree $T\in \mathcal{T}_{
\mathcal{M}}$ with $\mathcal{M}=\{1,2,3,4,5^3,6,7^2\}$ is illustrated in Figure \ref{T}.

For a sequence $\pi=\pi_1\pi_2\ldots\pi_n$,  the entry $\pi_i$, $1\leq i\leq n$, is said to be a {\em cyclic descent} (resp. a {\em cyclic  ascent}) if $\pi_i>\pi_{i+1}$ (resp. $\pi_i<\pi_{i+1}$) with the convention $\pi_{n+1}=\pi_1$. An entry $\pi_i$, $1\leq i\leq n$,  is said to be a {\em double cyclic descent} (resp. a {\em double cyclic  ascent}, a {\em cyclic peak}, a {\em cyclic valley} ) if  $\pi_{i-1}>\pi_i>\pi_{i+1}$ (resp. $\pi_{i-1}<\pi_i<\pi_{i+1}$,  $\pi_{i-1}<\pi_i>\pi_{i+1}$, $\pi_{i-1}>\pi_i<\pi_{i+1}$ ) with the convention $\pi_{n+1}=\pi_1$ and $\pi_{0}=\pi_{n}$. Let $cdes(\pi)$ (resp. $casc(\pi)$, $dcdes(\pi)$, $dcasc(\pi)$, $cpeak(\pi)$, $cval(\pi)$ ) denote the number of cyclic descents (resp. cyclic ascents, double cyclic descents, double cyclic ascents, cyclic peaks, cyclic valleys) of $\pi$.
Denote by $CDES(\pi)$ the set of cyclic descents of $\pi$.   For any sequence $\pi=\pi_1\pi_2\ldots \pi_n$, if $\pi_i$ is a   double cyclic  ascent or a   double cyclic descent,  the {\em cyclic $\pi_i$-factorization} of $\pi$ is defined to be $W_1\pi_iW_2W_3$, where
  $W_1$ (resp. $W_2$) is the maximal contiguous subsequence of $\pi$ immediately to the left (resp. right) of  $\pi_i$ whose entries  are all smaller than $\pi_i$ in the sense that the entries of $\pi$ are  arranged in a cycle clockwisely. For instance, let $\pi=15324$, then the cyclic $3$-factorization of $\pi$ is given by $W_13 W_2W_3$ where  $W_1$ is empty, $W_2=2$ and $W_3=415$.

In an ordered  labeled tree $T$, let $u$ be a  vertex of $T$ which is labeled by $v_0$. Suppose that the children of  $u$   are labeled by  $v_1, v_2, \ldots, v_{\ell}$ from left to right.  The number of  cyclic descents (resp. cyclic ascents) of $u$, denoted by $cdes(u)$ (resp. $casc(u)$),  is defined to be  the number  $cdes(v_0v_1v_2\ldots v_{\ell})$ (resp. $casc(v_0v_1v_2\ldots v_{\ell})$). The number of cyclic descents and  cyclic ascents  of $T$  are defined to be
$$
  cdes(T)=\sum_{u\in V(T)}cdes(u)$$
  and
  $$ casc(T)=\sum_{u\in V(T)}casc(u)
$$
 where $V(T)$ denotes the vertex set of $T$.
Denote by $eleaf(T)$ the number of leaves at even level.

Let $v$ be an odd vertex of a tree $T\in \mathcal{T}_n$ and let $u$ be its parent   labeled by $v_0$. Suppose that the children of the vertex $u$   are  labeled by $v_1, v_2, \ldots, v_{\ell}$  from left to right. Assume that the vertex $v$ is labeled by  $v_i$ for some $1\leq i\leq \ell$.  Then the odd vertex $v$ is said to be a {\em double cyclic descent } (resp. a {\em double cyclic ascent }, a {\em  cyclic peak }, a {\em  cyclic valley}) of $T$ if and only if  $v_i$ is a double cyclic descent (resp. a  double cyclic ascent, a  cyclic peak, a cyclic valley)  of the sequence $v_0v_1v_2\ldots v_\ell$.

Let $u$ be an even vertex of a tree $T\in \mathcal{T}_n$ which is labeled by $v_0$. Suppose that the children of $u$ are labeled by  $v_1, v_2, \ldots, v_{\ell}$   from left to right. Then the even vertex $u$ is said to be a {\em double cyclic descent } (resp. a {\em double cyclic ascent }, a {\em  cyclic peak},  a {\em  cyclic valley}) of $T$ if and only if  $v_0$ is a double cyclic descent (resp. a  double cyclic ascent, a  cyclic peak, a cyclic valley)  of the sequence $v_0v_1v_2\ldots v_\ell$.
Denote by $dcdes(T)$ (resp. $dcasc(T)$, $cpeak(T)$, $cval(T)$) the   number of
  double cyclic descents (resp. double cyclic ascents, cyclic peaks, cyclic valleys) of $T$.

 For example, if we let $T$ be a tree as shown  in Figure \ref{T}, we have $casc(T)=5$, $cdes(T)=4$,   $eleaf(T)=2$, $dcdes(T)=0$, $dcasc(T)=1$, $cpeak(T)=4$ and  $cval(T)=4$.

Now we are ready to state our main results of this section.

 \begin{theorem}\label{thtree}
   Let $$T_{\mathcal{M}}(x,y,z)=\sum_{T\in \mathcal{T}_{
 \mathcal{M}}}x^{casc(T)}y^{cdes(T)} z^{eleaf(T)},$$ where $\mathcal{M}=\{1^{k_1}, 2^{k_2}, \ldots, n^{k_n}\}$  and  $K=k_1+k_2+\ldots+k_n$ with $k_i\geq 1$. The polynomial
$T_{\mathcal{M}}(x,y,z)$ is  partial $\gamma$-positive and has the expansion

 $$
 T_{\mathcal{M}}(x,y,z)=\sum_{i=0}^{K-n}z^i \sum_{j=1}^{\lfloor{K+1-i\over 2}\rfloor} \gamma^{\mathcal{T}}_{\mathcal{M}, i,j}(xy)^j(x+y)^{K+1-i-2j},
 $$
 where
 $\gamma^{\mathcal{T}}_{\mathcal{M}, i,j}=|\{T\in \mathcal{T}_{\mathcal{M}}\mid cdes(T)=j, eleaf(T)=i, dcdes(T)=0\}|$.
\end{theorem}

\begin{figure}[h]
\centerline{\includegraphics[width=5cm]{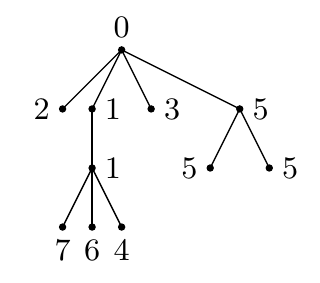}}
    \caption{A  tree $\mathcal{T}\in \mathcal{IT}_\mathcal{M}$ with $\mathcal{M}=\{1^2,2,3,4,5^3,6,7\}$.}\label{IT}
    \end{figure}

A tree $T\in \mathcal{T}_{\mathcal{M}}$ is said to be  {\em weakly increasing }  if  the labels of the vertices on the path from the root to any leaf is weakly increasing from top to bottom. See Figure \ref{IT} for an example. Denote by $\mathcal{IT}_{\mathcal{M}}$ the set of weakly increasing ordered labeled trees $T\in \mathcal{T}_{\mathcal{M}}$.

\begin{theorem}\label{thitree}
   Let $$IT_{\mathcal{M}}(x,y,z)=\sum_{T\in \mathcal{IT}_{
 \mathcal{M}}}x^{casc(T)}y^{cdes(T)} z^{eleaf(T)},$$ where $\mathcal{M}=\{1^{k_1}, 2^{k_2}, \ldots, n^{k_n}\}$  and  $K=k_1+k_2+\ldots+k_n$ with $k_i\geq 1$. The polynomial
$IT_{\mathcal{M}}(x,y,z)$ is  partial $\gamma$-positive and has the expansion

 $$
 IT_{\mathcal{M}}(x,y,z)=\sum_{i=0}^{K-n}z^i \sum_{j=1}^{\lfloor{K+1-i\over 2}\rfloor} \gamma^{\mathcal{IT}}_{\mathcal{M}, i,j}(xy)^j(x+y)^{K+1-i-2j},
 $$
 where
 $\gamma^{\mathcal{IT}}_{\mathcal{M}, i,j}=|\{T\in \mathcal{IT}_{\mathcal{M}}\mid cdes(T)=j, eleaf(T)=i, dcdes(T)=0\}|$.
\end{theorem}

In the following, we shall define  a group action  on ordered labeled trees in the sprit of  the Foata-Strehl
action  on permutations \cite{Fo}.

Let $T\in \mathcal{T}_{\mathcal{M}} $. For any vertex $u$ of $T$,
we define the   the {\em tree FS-action} $\psi_u$ on the  ordered labeled trees $T$ in the following way.

\begin{figure}[h]
\centerline{\includegraphics[width=10cm]{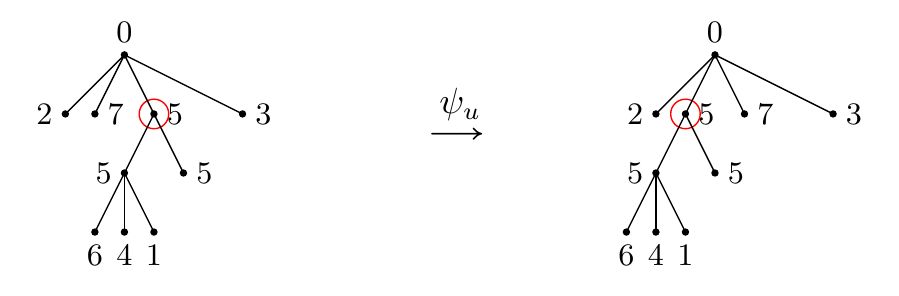}}
        \caption{An example of the action $\psi_u$ where the vertex $u$ is circled.}\label{F1}
    \end{figure}

\noindent{\bf Case 1.} If the vertex $u$ is an even leaf, a cyclic valley or a cyclic peak, let $\psi_u(T)=T$.
\noindent{\bf Case 2.} The vertex $u$ is a  double cyclic ascent or  a  double cyclic  descent.\\
\noindent{\bf Subcase 2.1.} The vertex $u$    is odd. Suppose that the   vertex $v$ labeled by $v_0$ is  the parent of $u$   and the children of $v$  are labeled by $v_1, v_2, \ldots, v_\ell$  from left to right. Assume that the vertex $u$ is labeled by $v_{k}$ for some $1\leq k\leq \ell $.  Let $T_i$, $(1\leq i\leq \ell)$, be the subtree rooted at the odd vertex labeled by $v_i$. Suppose that  the cyclic $v_k$-factorization of the sequence $v_0v_1v_2 \ldots v_\ell$ is given by $W_1 v_k W_2W_3$. Then let $\psi_u(T)$ be the tree obtained from $T$ by rearranging the subtrees $T_1, T_2, \ldots, T_{\ell}$ such that   $CDES(W_2v_kW_1W_3)=CDES(v_0v'_1v'_2\ldots v'_\ell)$, where  the children of the  vertex $v$ in $\psi_u(T)$ are labeled by $v'_1, v'_2, \ldots, v'_\ell$  from left to right. See Figure \ref{F1} for an example.

\begin{figure}[h]
\centerline{\includegraphics[width=10cm]{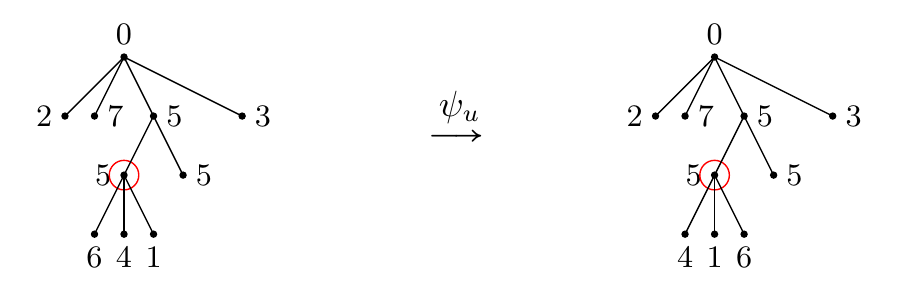}}
        \caption{An example of the action $\psi_u$ where the vertex $u$ is circled.}\label{F2}
    \end{figure}

  \noindent{\bf Subcase 2.2.} The vertex $u$ is    even. Suppose that the vertex $u$ is labeled by $v_0$ and the   children of $u$ are  labeled by $v_1, v_2, \ldots, v_\ell$   from left to right.  Let $T_i$, $(1\leq i\leq \ell)$, be the subtree rooted at the odd vertex labeled by $v_i$. Suppose that  the cyclic $v_0$-factorization of the sequence $v_0v_1v_2 \ldots v_\ell$ is given by $W_1 v_0 W_2W_3$. Then let $\psi_u(T)$ be the tree obtained from $T$ by rearranging the subtrees $T_1, T_2, \ldots, T_{\ell}$ such that    $CDES(W_2v_0W_1W_3)=CDES(v_0v'_1v'_2\ldots v'_\ell)$, where the children of the  vertex $u$ in $\psi_u(T)$ are labeled by $v'_1, v'_2, \ldots, v'_\ell$  from left to right. See Figure \ref{F2} for an example.

The following observation will play an essential role in the proof of Theorem \ref{thtree}.
\begin{observation}
Let $\mathcal{M}=\{1^{k_1}, 2^{k_2}, \ldots, n^{k_n}\}$  and  $K=k_1+k_2+\ldots+k_n$ with $k_i\geq 1$.
Then we have
\begin{equation}\label{eq1}
eleaf(T)+dcdes(T)+dcasc(T)+cpeak(T)+cval(T)=K+1,
\end{equation}
   \begin{equation}\label{eq2}
 cpeak(T)=cval(T),
\end{equation}
\begin{equation}\label{eq3}
 cdes(T)=cval(T)+dcdes(T)
\end{equation}
and
\begin{equation}\label{eq4}
 casc(T)=cpeak(T)+dcasc(T)
\end{equation}
for any $T\in \mathcal{T}_{\mathcal{M}}$.

\end{observation}
 We are now ready for the proof of Theorem \ref{thtree}.

 \noindent{\bf Proof of Theorem \ref{thtree}.} Define $\mathcal{T}_{\mathcal{M}, i}$ be the set of ordered labeled trees $T\in \mathcal{T}_{\mathcal{M}}$ with $eleaf(T)=i$.
 Then Theorem \ref{thtree} is equivalent to
 \begin{equation}\label{eqT}
 \begin{array}{lll}
 \sum_{T\in \mathcal{T}_{\mathcal{M}, i}}x^{casc(T)}y^{cdes(T)}&=&  \sum_{T\in \widetilde{ \mathcal{T}_{\mathcal{M}, i}}}(xy)^{cdes(T)}(x+y)^{K+1-i-2cdes(T)}.
 \end{array}
 \end{equation}
 where $\widetilde{ \mathcal{T}_{\mathcal{M}, i}}$ denotes the set of ordered labeled trees $T\in \mathcal{T}_{\mathcal{M}, i}$ with $dcdes(T)=0$.

 Clearly, the tree FS-actions $\psi_u$'s are involutions and  commute.  Thus, for any  $S\subseteq V(T)$, we can define the function $\psi_{S}: \mathcal{T}_{\mathcal{M}} \rightarrow \mathcal{T}_{\mathcal{M}} $ by $\psi_{S}=\prod_{u\in S}\psi_u$, where the product is the functional compositions.  Hence the group  $\mathbb{Z}_2^{K+1}$ acts on  $\mathcal{T}_{\mathcal{M}}$ via the function $\psi_S$. Since
the statistic  ``eleaf" is invariant under this group action, it divides $\mathcal{T}_{\mathcal{M}, i}$ into some disjoint
orbits. For each $T\in \mathcal{T}_{\mathcal{M}, i} $,  let $Orb(T) = \{g(T) : T \in  \mathbb{Z}_2^{K+1}\}$ be the orbit of $T$ under the
tree FS-action. Notice that the vertex $u$ is a double cyclic descent of $T$ if and only if the vertex $u$ is a double cyclic ascent of $\psi_u(T)$.  This yields that there is a unique tree $\widetilde{T}\in Orb(T) $ such that $dcdes(T)=0$. Therefore, we have
$$
\begin{array}{lll}
\sum_{T\in   Orb(\widetilde{T})}x^{casc(T)}y^{cdes(T)} &=& x^{cpeak(\widetilde{T})}y^{cval(\widetilde{T})}(x+y)^{cdasc(\widetilde{T})}\\
&=& x^{cpeak(\widetilde{T})}y^{cval(\widetilde{T})}(x+y)^{K+1-i-cpeak(\widetilde{T})-cval(\widetilde{T})}\\
 &=& (xy)^{cdes(\widetilde{T})}(x+y)^{K+1-i-2cdes(\widetilde{T})},\\
\end{array}
$$
where the second equality follows from (\ref{eq1}) and the third equality follows from (\ref{eq2})-(\ref{eq4}).  Summing over all orbits of $\mathcal{T}_{\mathcal{M}, i}$ under the tree FS-action then gives (\ref{eqT}), completing the proof. \qed

\noindent{\bf Proof of Theorem \ref{thitree}.} From the construction of the cyclic FS-action $\psi_u$,   it is easily seen that $\psi_u(T)\in \mathcal{Q}_{\mathcal{M}} $ for  any $T\in \mathcal{IT}_{\mathcal{M}}$. Therefore, Theorem \ref{thitree} can be proved  by similar reasoning as in the proof of Theorem \ref{thtree}. \qed
\section{Partial $\gamma$-positivity for  Quasi-Stirling Permutations }
Let $\pi$ be a permutation, the leftmost entry of $\pi$ is denoted by $first(\pi)$. Similarly, for an ordered labeled tree $T$, we denote by $first(T)$  the label of the leftmost   child of the root.

In \cite{Yanzhu},  Yan and Zhu   established a bijection $\phi$  between $\mathcal{T}_{\mathcal{M}}$ and $
\overline{\mathcal{Q}}_{\mathcal{M}}$. First we give an overview of the bijection $\phi$.
If $T$ has only one vertex, let $\phi(T)=\epsilon$, where $\epsilon$ denotes the empty permutation.  Otherwise, suppose that $first(T)=r$.

\noindent {\bf Case 1.}  The leftmost child of the root is a leaf. Let $T_0$ be the tree obtained from $T$ by removing  the leftmost child of the root together with  the edge incident to it. Define
    $\phi(T)=r\phi(T_0)$.

    \noindent {\bf Case 2.} The leftmost child of the root has $k$ children. For $1\leq i\leq k$,  let $T_i$ be the subtree rooted at the $i$-th child of the leftmost child of the root.
         Denote by $T'_i$ the tree obtained from $T_i$ by relabeling its root by $0$.
         Let $T_0$ be the tree obtained from $T$ by removing all the subtrees $T_1$, $T_2$, $\ldots,$ $T_k$ and  all the vertices labeled by $r$ together with  the edges incident to them. Define $\phi(T)=r\phi(T'_1)r\phi(T'_2)\ldots r\phi(T'_k)r\ \phi(T_0)$.

For instance, let $T$ be a tree illustrated in Figure \ref{T}. By applying the map  $\phi$, we have $\phi(T)=2773516455$.

In \cite{Yanyang}, Yan et al. proved that the bijection $\phi$ has the following property.

\begin{lemma}{\upshape   ( \cite{Yanyang}) }\label{lemphi}
For any nonempty multiset $\mathcal{M}$, the map $\phi$ induces a bijection between $\mathcal{T}_{\mathcal{M}}$ and $
\overline{\mathcal{Q}}_{\mathcal{M}}$  such that $$(cdes, casc, eleaf)T= (des, asc, plat)\phi(T)$$ for any $T\in \mathcal{T}_{\mathcal{M}}$.
\end{lemma}

Combining Lemma \ref{lemphi} and Theorem \ref{thtree}, we immediately derive the following partial $\gamma$-positivity for $\overline{Q}_{\mathcal{M}}(x,y,z)$, thereby confirming Conjecture 
\ref{con1}. 
\begin{theorem}\label{thquasi}
Let $\mathcal{M}=\{1^{k_1}, 2^{k_2}, \ldots, n^{k_n}\}$  and  $K=k_1+k_2+\ldots+k_n$ with $k_i\geq 1$.
The polynomial
$\overline{Q}_{\mathcal{M}}(x,y,z)$ is  partial $\gamma$-positive and has the expansion
 $$
 \overline{Q}_{\mathcal{M}}(x,y,z)=\sum_{i=0}^{K-n}z^i \sum_{j=1}^{\lfloor{K+1-i\over 2}\rfloor} \gamma^{\mathcal{T}}_{\mathcal{M}, i,j}(xy)^j(x+y)^{K+1-i-2j},
 $$
 where
 $\gamma^{\mathcal{T}}_{\mathcal{M}, i,j}=|\{T\in \mathcal{T}_{\mathcal{M}}\mid cdes(T)=j, eleaf(T)=i, dcdes(T)=0\}|$.
\end{theorem}

From the construction of the map $\phi$, it is easily seen that $\phi(T)\in \mathcal{Q}_{\mathcal{M}} $ for any $T\in \mathcal{IT}_{\mathcal{M}}$. Hence, the following result follows immediately from Lemma \ref{lemphi}.
\begin{lemma}\label{lemphi1}
For any nonempty multiset $\mathcal{M}$, the map $\phi$ induces a bijection between $\mathcal{IT}_{\mathcal{M}}$ and $
 \mathcal{Q}_{\mathcal{M}}$  such that $$(cdes, casc, eleaf)T= (des, asc, plat)\phi(T)$$ for any $T\in \mathcal{IT}_{\mathcal{M}}$.
\end{lemma}

Combining   Lemma  \ref{lemphi1} and Theorem \ref{thitree},  we are led to an alternative proof of the 
partial $\gamma$-positivity of the polynomial $Q_{\mathcal{M}}(x,y,z)$. 

\begin{theorem}\label{ths}
Let $\mathcal{M}=\{1^{k_1}, 2^{k_2}, \ldots, n^{k_n}\}$  and  $K=k_1+k_2+\ldots+k_n$ with $k_i\geq 1$.
The polynomial
$Q_{\mathcal{M}}(x,y,z)$ is  partial $\gamma$-positive and has the expansion
 $$
 Q_{\mathcal{M}}(x,y,z)=\sum_{i=0}^{K-n}z^i \sum_{j=1}^{\lfloor{K+1-i\over 2}\rfloor} \gamma^{\mathcal{IT}}_{\mathcal{M}, i,j}(xy)^j(x+y)^{K+1-i-2j},
 $$
 where
 $\gamma^{\mathcal{IT}}_{\mathcal{M}, i,j}=|\{T\in \mathcal{IT}_{\mathcal{M}}\mid cdes(T)=j, eleaf(T)=i, dcdes(T)=0\}|$.
\end{theorem}
 
\begin{remark}
In \cite{Lin}, Lin, Ma and Zhang provided three combinatorial interpretations for the $\gamma$-coefficients of the partial $\gamma$-positivity expansion of $Q_{\mathcal{M}}(x,y,z)$.  Theorem  \ref{ths} gives a new combinatorial interpretation.
\end{remark}

We conclude this section with an alternative proof of the partial $\gamma$-positivity of the polynomial $\overline{Q}_{\mathcal{M}}(x,y,z)$ by employing the result obtained by Yan  et al. \cite{Yanyang} and the
partial $\gamma$-positivity of the polynomial $ Q_{\mathcal{M}}(x,y,z)$ proved by Lin, Ma and Zhang \cite{Lin}.
In \cite{Yanyang},   Yan  et al.   obtained the following result concerning the polynomial $\overline{Q}_{\mathcal{M}}(x,y,z)$.

\begin{theorem}{\upshape   ( \cite{Yanyang}) }\label{thyan}
Let $\mathcal{M}=\{1^{k_1}, 2^{k_2}, \ldots, n^{k_n}\}$, $\mathcal{M'}=\{1^{K-n+1}, 2,3,\ldots, n\}$ and $K=k_1+k_2+\ldots+k_n$ with $k_i\geq 1$. We have
      $\overline{Q}_{\mathcal{M}}(x,y,z)=\overline{Q}_{\mathcal{M'}}(x,y,z)$.
    \end{theorem}
Notice that quasi-Stirling permutations of $\mathcal{M'}=\{1^{K-n+1}, 2,3,\ldots, n\}$ are equivalent to Stirling permutations of $\mathcal{M'}=\{1^{K-n+1}, 2,3,\ldots, n\}$. Thus, by Theorem \ref{thyan}, we have $\overline{Q}_{\mathcal{M}}(x,y,z)=Q_{\mathcal{M'}}(x,y,z)$.  Thus, the partial $\gamma$-positivity of the polynomial $ Q_{\mathcal{M}}(x,y,z)$ implies the the partial $\gamma$-positivity of the polynomial $ \overline{Q}_{\mathcal{M}}(x,y,z)$.

 \noindent{\bf Acknowledgments.}
 This work was supported by  the National Natural Science Foundation of China (12071440).


\end{document}